\documentclass[12pt]{amsart}
\usepackage[english]{babel}
\usepackage{amsmath}
\usepackage{amsthm}
\usepackage{amssymb}
\usepackage{mathrsfs}
\usepackage{mathdots}
\usepackage{enumerate,MnSymbol}
\usepackage[notcite, final, notref]{showkeys}
\usepackage{dsfont}
\usepackage{color}
\usepackage{shadow}
\newtheorem{theorem}{Theorem}[section]

\newtheorem{lemma}[theorem]{Lemma}
\newtheorem{proposition}[theorem]{Proposition}

\newtheorem{definition}[theorem]{Definition}

\theoremstyle{definition}
\newtheorem{remark}[theorem]{Remark}

\newcommand{\ov}{\overline}

\newcommand{\inner}[1]{\left\langle  #1 \right\rangle }

\newcommand{\fact}[1]{\left[  #1 \right]_q ! \, }

\title[Generalized $q$-Fock spaces and structural identities]{Generalized $q$-Fock spaces and structural identities}

\author[D. Alpay]{Daniel Alpay}
\address{(DA)
Faculty of Mathematics, Physics, and Computation\\
Schmid College of Science and Technology\\
Chapman University\\
One University Drive
Orange, California 92866\\
USA}
\email{alpay@chapman.edu}

\author[P. Cerejeiras]{Paula Cerejeiras}
\address{(PC) Center for research and development in mathematics and applications\\Department of mathematics, University of Aveiro  \\
Campus Universit\'ario de Santiago  \\ 3810-193 Aveiro\\Portugal}
\email{pceres@ua.pt}

\author[U. Kaehler]{Uwe Kaehler}
\address{(UK) Center for research and development in mathematics and applications\\Department of mathematics, University of Aveiro  \\
Campus Universit\'ario de Santiago  \\ 3810-193 Aveiro\\Portugal}
\email{ukaehler@ua.pt}

\author[B. Schneider]{Baruch Schneider}
\address{(BS) University of Ostrava\\ Department of Mathematics\\
30.dubna 22, 70200 Ostrava\\Czech Republic }
\email{baruch.schneider@osu.cz}

\begin{document}

\begin{abstract}
  Using $q$-calculus we study a family of reproducing kernel Hilbert spaces which interpolate between the Hardy space and the Fock space. We give characterizations of these spaces in terms of classical operators such as integration and backward-shift operators, and their $q$-calculus counterparts. Furthermore, these new spaces allow us to study intertwining operators between classic backward-shift operators and the q-Jackson derivative.
\end{abstract}

\thanks{D. Alpay thanks the Foster G. and Mary McGaw Professorship in
  Mathematical Sciences, which supported his research.}

\maketitle

\section{Introduction}
\setcounter{equation}{0}
\subsection{Prologue}
The Hardy space of the open unit disk $\mathbb D$, here denoted by $\mathbf H_2 = \mathbf H_2(\mathbb D)$, is the reproducing kernel Hilbert space with reproducing kernel
\[
  \frac{1}{1-z\overline{w}}=\sum_{n=0}^\infty z^n\overline{w}^n, \quad z,w\in\mathbb D,
\]
and plays a key role in operator theory, linear system theory and Schur analysis. On the other hand, the Bargmann-Fock-Segal space, here
denoted by $\mathcal F$ and called Fock space for short, is the reproducing kernel Hilbert space with reproducing kernel
  \[
e^{z\overline{w}}=\sum_{n=0}^\infty \frac{z^n\overline{w}^n}{n!},\quad z,w\in\mathbb C,
    \]
    and plays a key role in quantum mechanics (and more recently in signal processing).\\

    The Hardy space $\mathbf H_2$ can be characterized (up to a positive multiplicative factor for the inner product)
    as the only Hilbert space of power series converging at the origin and such that
\begin{equation}
  R_0^*=M_z,
\end{equation}
  where $M_z$ is the operator of multiplication by $z$ and
\begin{equation}
  R_0f(z)=\frac{f(z)-f(0)}{z}.
  \end{equation}
  Note that in $\mathbf H_2$ we have the identities
  \begin{equation}
    R_0R_0^*= \mathcal{I}, \quad \mbox{\rm and }\quad
R_0M_z-M_zR_0=\mathcal{I}-R_0^*R_0=C^*C,
\end{equation}
where $Cf=f(0)$ and $\mathcal{I}$ is the identity operator. We remark that
\begin{equation}
 \mathcal{I}-R_0^*R_0=C^*C,
\end{equation}
  which we will call  {\sl structural identity}, is the simplest of a family of identities characterizing de Branges spaces.\smallskip

  Similarly, and besides Bargmann celebrated characterization $\partial^*=M_z$ (see \cite{MR0157250,bargmann}), the Fock space is (still up to a
  positive multiplicative factor for the inner product) the only Hilbert space of power series converging at the origin and such that
\begin{equation}
  R_0^*={\rm I},
\end{equation}
where ${\rm I}$ is the integration operator (see \cite{Trevor23})
\begin{equation}
({\rm I}f)(z)=  \int_{[0,z]}f(s)ds.
  \end{equation}

  \subsection{The paper}
    The $q$-calculus allows to define a continuum of spaces between
    $\mathbf H_2$ and $\mathcal F$, namely the family of reproducing kernel Hilbert spaces $\mathbf H_{2,q}$ indexed by $q\in[0,1]$ and
    with reproducing kernel
    \[
K_q(z,w)=\sum_{n=0}^\infty \frac{z^n\overline{w}^n}{[n]_q!},\quad q\in[0,1],\quad z,w\in\mathbb D_{1/1-q},
      \]
where in the above expression
      \[
       \mathbb D_{1/1-q}=\begin{cases}\,  \mathbb D_{\infty}= \mathbb C,\quad\hspace{1cm} q=1,\\
                              \,\left\{z\in\mathbb C\,:\, |z|<\frac{1}{1-q}\right\},\quad q\in[0,1)\end{cases}.
        \]
Furthermore, $\fact{0}=1$ and $
[n]_q!=1\cdot(1+q)\cdot(1+q+q^2)\cdots (1+q+\cdots +q^{n-1}),\quad n\in\mathbb N.$ Thus, in this notation, we have
        \[
        \mathbf H_{2,0}=\mathbf H_2\quad{\rm and}\quad \mathbf H_{2,1}=\mathcal F,
      \]
with
      \[
          K_{0}(z,w):= k_{2,0}(z,w)=\frac{1}{1-z\overline{w}}\quad {\rm and}\quad  K_{1}(z,w):= k_{2,1}(z,w)=         e^{z\overline{w}}.
          \]

        The $q$-calculus allows to gather into a common umbrella problems pertaining to the classical Hardy space $\mathbf H_2$ of the open unit
disk and problems pertaining to the Fock space.

Consider now
\begin{equation}
  \label{rq}
R_qf(z)=\frac{f(z)-f(qz)}{(1-q)z}, \qquad 0\leq q < 1,
\end{equation}
while for $q=1$, we consider $R_1=\partial.$ In this way we have a progression between two fundamental linear operators in analysis, namely the backward-shift and the differentiation operators.
  Then, one can introduce the $q$-Fock space $\mathbf H_{2,q}$ as the unique (up to a multiplicative positive constant) space of power series such that $R_q^*=M_z$. The case $q=1$ corresponds to the classical Fock space (see  \cite{MR0157250}). It is important to note already at this stage that these operators satisfy a $q$-commutator relation (see also Lemma \ref{lemma-q})
  \begin{equation}
    R_qM_z-qM_zR_q= \mathcal{I}.
\end{equation}

%
%

 \section{$q$-calculus}

 \subsection{Iterative powers of the operator $R_q$}


Recall that $R_q$ was defined by \eqref{rq}.

\begin{proposition} Let ${\Lambda_q} f(z) = f(qz)$. We have
\begin{equation}
R_q^nf(z)=\frac{\Pi_{k=1}^n (1-q^k {\Lambda_q})}{(1-q)^n} R_0^nf(z), \qquad 0\leq q < 1,\quad n=1,2,\ldots
  \end{equation}
\end{proposition}
\begin{proof} Firstly, we observe the intertwining between $R_0$ and $\Lambda_q,$ 
\begin{gather*}
R_0 {\Lambda_q}f(z) = R_0f(qz)= \frac{f(qz)-f(0)}{z} = q \frac{f(qz)-f(0)}{qz} = q {\Lambda_q} R_0 f(z).
\end{gather*}
Secondly,
\[
\begin{split}
  R_qf(z)& = \frac{f(z)-f(qz)}{(1-q)z} = \frac{f(z)-f(0) -f(qz) + f(0)}{(1-q)z} \\
&= \frac{1}{1-q}\frac{f(z)-f(0) }{z}  - \frac{q}{1-q} \frac{f(qz) - f(0)}{qz} \\
&= \frac{(1-q {\Lambda_q})R_0}{1-q} f(z).
\end{split}
\]
Hence,
\[
  \begin{split}
    R_q^2f(z) &=  \left(  \frac{(1-q {\Lambda_q})R_0}{1-q} \right)^2 f(z)\\
    &=  \frac{(1-q {\Lambda_q})R_0 (1-q{\Lambda_q})R_0}{(1-q)^2}f(z)\\
&= \frac{(1-q {\Lambda_q}) (1-q^2 {\Lambda_q})R_0^2}{(1-q)^2}f(z),\\
\end{split}
\]
and by induction the result holds:
\[
  \begin{split}
    R_q^nf(z)& =  \left(  \frac{(1-q {\Lambda_q})R_0}{1-q}\right)^n f(z)\\
    &=  \frac{(1-q {\Lambda_q})R_0(1-q {\Lambda_q})R_0 \cdots (1-q{\Lambda_q})R_0}{(1-q)^n}f(z)\\
&= \frac{(1-q {\Lambda_q}) (1-q^2 {\Lambda_q})R_0^2\cdots (1-q{\Lambda_q})R_0 }{(1-q)^n}f(z)\\& \hspace{2mm} \vdots\\
&= \frac{(1-q {\Lambda_q}) (1-q^2 {\Lambda_q}) \cdots (1-q^n{\Lambda_q})R^n_0 }{(1-q)^n}f(z).
\end{split}
\]
\end{proof}

As we will see later (see Theorems \ref{Th:4.3} and \ref{Th:4.4}), $R_q^*$ has completely different properties depending on which of the spaces at hand we compute the adjoint.

\subsection{$q-$Stirling numbers associated to higher commutation relations} In this subsection let us recall some facts regarding higher-commutator relations in $q$-calculus. While they can be found, e.g., in~\cite{Schork2016} for the sake of self-sufficiency of the paper we present  them with proofs. 

For the $q$-commutator we have the following well-known formula.

\begin{lemma}[$q$-commutator] For the $q$-commutator it holds the following identity:
\begin{equation}
[R_q, M_z]_q :=  R_q M_z - q M_z R_q = \mathcal{I}.
\end{equation}
\label{lemma-q}
\end{lemma}

\begin{proof} We have
\begin{gather*}
R_q M_z z^n  =  R_q z^{n+1} = (1+q+\cdots + q^n) z^{n}, \quad n=0,1,2, \ldots
\end{gather*} while
\begin{gather*}
q M_z R_q z^n = q M_z R_q 1 = 0,\quad n=0, \\
qM_z R_q z^n  =  q M_z (1+q+\cdots + q^{n-1}) z^{n-1} = (q+q+\cdots + q^n) z^{n}, \quad n=1,2, \ldots
\end{gather*} so that it holds $(R_q M_z - q M_z R_q) z^n =z^n,$ for all $n \in \mathbb{N}_0.$
\end{proof}

We define our $q-$Stirling numbers as coefficients $S(n,k)$ of the following commutation relation (see \cite{MR3816055}):
\begin{equation}\label{n-commutator}
(M_z R_q)^n := \sum_{k=1}^n S(n, k) M_z^k  R_q^k,\qquad n \in \mathbb{N}.
\end{equation}

This formula can also be found in \cite{Schork2016} (Theorem 3.1) and indirectly also in \cite{Milne1978}. Furthermore, in \cite{Schork2016}, Section 4.1. there is a general exposition on how to construct such higher order commutator relations including formulae for terms of the type $(M_z^r R_q^s)^n$ with $r, s$ multi-indices. 

\begin{lemma} We have for these $q-$Stirling numbers  the following recursion formula
\begin{gather*}
S(1,1) = 1; \\
S(n, n) = S(n-1, n-1) q^{n-1}, \qquad n=2, 3, \ldots;\\
S(n, k) = (1+q+\cdots +q^{k-1}) S(n-1, k) + q^{k-1} S(n-1, k-1), \qquad k=2, \ldots, n-1.
\end{gather*}
\end{lemma}

This recursion formula is known in the literature. One can find it in~\cite{Milne1978} formula (1.15) on page 93) or in the book~\cite{Schork2016} (Section 3.3, page 68 onwards).

\begin{proof} In order to simplify notation, we write the expression for the $q-$Stirling numbers as
$$(ab)^n := \sum_{k=1}^n S(n, k) a^k  b^k.$$
From the $q-$commutator we get $ba = 1+q ab$ so that
\begin{gather*}
b^n a = b^{n-1} (ba) = b^{n-1} (1+q ab) = b^{n-1} + q (b^{n-1}a) b \\
= b^{n-1} + q \left[ b^{n-2} + q (b^{n-2}a) b  \right] b = (1+q) b^{n-1} + q^2 (b^{n-2}a) b^2 \\
\vdots \\
=  (1+q+ \cdots + q^{n-1}) b^{n-1} + q^n a b^n.
\end{gather*}
Replacing in the above formula for the $q-$Stirling numbers we obtain
\begin{gather*}
(ab)^n = \sum_{k=1}^n S(n, k) a^k  b^k \\
= (ab)^{n-1}(ab) = \big[ \sum_{k=1}^{n-1} S(n-1, k) a^k  b^k \big] (ab) \\
= \sum_{k=1}^{n-1} S(n-1, k) a^k  (b^k a)b = \sum_{k=1}^{n-1} a^k S(n-1, k)   \left[  (1+q+ \cdots + q^{k-1}) b^{k-1} + q^k a b^k \right] b
\end{gather*}
\begin{gather*}
= \sum_{k=1}^{n-1}   \left[  (1+q+ \cdots + q^{k-1}) S(n-1, k)  a^k  b^{k} + q^k S(n-1, k)  a^{k+1} b^{k+1} \right]\\
= \sum_{k=1}^{n-1}  (1+q+ \cdots + q^{k-1}) S(n-1, k)  a^k  b^{k} + \sum_{k=2}^{n}   q^{k-1} S(n-1, k-1)  a^{k} b^{k} \\
= S(n-1, 1) ab + \sum_{k=2}^{n-1} \left[ (1+q+ \cdots + q^{k-1}) S(n-1, k) + q^{k-1} S(n-1, k-1)    \right] a^k b^k + q^{n-1} S(n-1, n-1)  a^{n} b^{n},
\end{gather*}so that we have $S(1,1)=1,$
$$S(n,1) = S(n-1, 1), \quad S(n, n) = q^{n-1} S(n-1, n-1),$$
for $n=2, 3, \ldots$ and
$$S(n,k) = (1+q+ \cdots + q^{k-1}) S(n-1, k) + q^{k-1} S(n-1, k-1),$$
for $k=2, \ldots, n-1.$
\end{proof}

One can easily see the first $q-$Stirling numbers

\begin{center}
\begin{tabular}{|c||c|c|c|c|c|}
\hline
$S(n,k)$ & 1& 2& 3& 4 \\
\hline
\hline
1 & 1& & &  \\
\hline
2& 1& $q$ & &  \\
\hline
3& 1& $2q+q^2$ & $q^3$ & \\
\hline
4& 1& $q^3+3q^2 +3q$ & $q^5+2q^4+3q^3$ & $q^6$ \\
\hline
\end{tabular}
\end{center}

\begin{remark}
  We need to point out that there are two types of $q$-Stirling numbers of the first or of the second kind in the literature. The more classic ones were obtained by studying the corresponding partition problems in $q$-calculus~(for a review on this topic see \cite{Cai2, Cai1}). Here we have them as coefficients of the expansion of $(M_z R_q)^n$ in (\ref{n-commutator}) in the same way as in~\cite{Milne1978} and~\cite{Schork2016}.
 Only in the classic case of $q=1$ this type of coefficients coincides with classic Stirling numbers of the second kind, i.e. with the numbers of partitions of a set of $n$ objects into $k$ non-empty subsets.
\end{remark}

\section{The $q$-Fock space}
\setcounter{equation}{0}
Consider the positive definite function $E_q(z\overline{w})$ given by the $q$-exponential:
\begin{equation}
  E_q(z)=\sum_{k=0}^\infty \frac{z^k }{\fact{k}}=\frac{1}{\prod_{j=0}^\infty (1-z(1-q)q^j)}=: \frac{1}{(z(1-q);q)_\infty}, \quad z  \in \mathbb{D}_{1/1-q},
  \label{groningen}
\end{equation} evaluated at $z\overline w,$
with $[0]_q=1$ and $[k]_q=1+q+\cdots +q^{k-1}$ for $k=1,2, \ldots,$ and $\fact{k}=\prod_{j=0}^k[j]_q$, i.e.
\[
\fact{k}= [1]_q [2]_q \cdots [k]_ q = 1\cdot (1+q)\cdot (1+q+q^2)\cdots (1+q+\cdots +q^{k-1}).
\]
The term $(a;q)_n=\prod_{j=0}^{n-1}(1-aq^j)$ denotes the $q$-Pochhammer symbol.

\begin{definition}
  We denote by $\mathbf H_{2,q}$ the reproducing kernel Hilbert space of functions analytic in $|z|<\frac{1}{1-q}$ with reproducing kernel $E_q(z\overline{w})$.
  \end{definition}

As stated before when $q=0$ we get back the classical Hardy space of the open unit disk, while $q\rightarrow 1$ leads to the classical Fock space; see e.g. \cite{duren,MR1102893,MR924157} for the former, \cite{Zhu2012} for the latter.


For functions belonging to the $q$-Fock space we have the following characterization based on its power series expansion.
\begin{lemma}
  $f(z)=\sum_{n=0}^\infty a_nz^n$ belongs to   $\mathbf H_{2,q}$ if and only if
  \begin{equation}
\sum_{n=0}^\infty \fact{n} \, |a_n|^2<\infty.
    \end{equation}
  \end{lemma}

Based on the $q$-Jackson integral (see \cite{Jackson1910}, \cite{Kac2001})
$$
\int_0^a f(x)d_qx:=(1-q)a\sum_{k=0}^\infty q^k f(q^k a),
$$
we can define the following $q$-integral transform. 

\begin{definition}\label{Def:3.01}
Given a bounded function $f:[0,-1+ 1/(1-q)] \to \mathbb{R}$ we define its $q$-integral transform as
$$
\mathcal{M}_q f(z)=\int_0^{1/(1-q)} t^{z-1}f(qt)d_qt := \sum_{k=0}^\infty q^k \Big( \frac{q^k}{1-q} \Big)^{z-1} f \Big( \frac{q^{k+1}}{1-q} \Big).
$$
\end{definition}


With the help of this $q$-integral transform we get that the coefficients $\frac{1}{\fact{n}}$ satisfy the moment problem
%
\begin{eqnarray*}
\fact{n} = \mathcal{M}_q(E_q^{-1})(n+1) & = & \int_0^{1/(1-q)} t^{n}E_q^{-1}(qt)d_qt \\
& =&  \frac{(q;q)_\infty}{(1-q)^n} \sum_{k=0}^\infty \frac{q^{(n+1)k}}{(q;q)_k}, \quad 0< q <1, ~ n \in \mathbb{N}_0,
\end{eqnarray*}
where  $(a;q)_n = \Pi_{j=0}^{n-1} (1-aq^j)$ denotes the $q$-Pochhammer symbol and $(q;q)_n=\frac{(q;q)_\infty}{(q^{n+1};q)_\infty}.$

For the disk $\mathbb{D}_{1/(1-q)}$ we have the measure (see \cite{Leeuwen1995})
$$d\mu_q(z) = (q;q)_\infty \sum_{k=0}^\infty \frac{q^{k}}{(q;q)_k} d\lambda_{r_k}(z),$$
where $r_k = \frac{q^{k/2}}{\sqrt{1-q}}$ while $d\lambda_{r_k}$ is the normalized Lebesgue measure in the circle of radius $r_k.$ This leads to the following characterization of the space  $\mathbf H_{2,q}$.


\begin{theorem}
The space $\mathbf H_{2,q}$ corresponds to the space of all analytic functions in the disk $\mathbb D_{\frac{1}{1-q}}=\{z :  |z|<\frac{1}{1-q}\}$ satisfying the condition
$$
\iint_{\mathbb D_{\frac{1}{1-q}}} |f(z)|^2 d\mu_q(z)<\infty. $$
The inner product of $\mathbf H_{2,q}$ is given by
$$
\frac{1}{2\pi} \iint_{\mathbb D_{\frac{1}{1-q}}}  f(z)\overline{g(z)}  d\mu_q(z) =\sum_{n=0}^\infty f_n \overline{g_n}\fact{n}.
$$
\end{theorem}

\begin{proof} We have
\begin{eqnarray*}
\frac{1}{2\pi} \iint_{\mathbb D_{\frac{1}{1-q}}} z^n \overline{z}^m  d\mu_q(z) & = & \frac{ (q;q)_\infty}{2\pi} \sum_{k=0}^\infty \frac{q^{k}}{(q;q)_k} r_k^{n+m}
 \underbrace{\int_0^{2\pi} e^{i (n-m) \theta} d\theta}_{= 2\pi \delta_{n,m}} \\
 & = &  \delta_{n,m} (q;q)_\infty \sum_{k=0}^\infty \frac{q^{k}}{(q;q)_k}  r_k^{2n}
 \\
  & = & \delta_{n,m} \frac{(q;q)_\infty}{(1-q)^n} \sum_{k=0}^\infty \frac{q^{(n+1)k}}{(q;q)_k}.
\end{eqnarray*}
\end{proof}

Combining this result with our moment problem we obtain
\begin{gather*}
\fact{n} = \mathcal{M}_q(E_q^{-1})(n+1)
=\int_0^{1/(1-q)} t^{n}E_q^{-1}(qt)d_qt
= \frac{\delta_{n,m}}{2\pi} \iint_{\mathbb D_{\frac{1}{1-q}}}  z^n \overline z^m   d\mu_q(z).
\end{gather*}
We observe that for $q \to 1$ we obtain $d\mu_q(z) = \frac{1}{2} e^{-|z|^2} dx dy.$ \\

Also we get a convolution-type formula for our $q$-integral transform.

\begin{lemma}
Given bounded functions $f_1, f_2:[0,-1+ 1/(1-q)] \to \mathbb{R}$ it holds (pointwisely)
\begin{equation} \label{Eq:3.0001}
\mathcal M_q(f_1)(z) \mathcal M_q(f_2)(z) = \left(\frac{1}{1-q} \right)^{z-1} \mathcal M_q(f_1 \circ f_2)(z),
\end{equation} where
\begin{equation} \label{Eq:3.0001a}
f_1 \circ f_2 \left( q\frac{q^m}{1-q} \right) := \sum_{k=0}^m  f_1  \left( \frac{q^{k+1}}{1-q} \right)f_2  \left( \frac{q^{m+1-k}}{1-q} \right).
\end{equation}
\end{lemma}

\begin{proof} From Definition \ref{Def:3.01} we have
\begin{gather*}
\mathcal M_q(f_1)(z) \mathcal M_q(f_2)(z) = \left( \sum_{k=0}^\infty q^k \Big( \frac{q^k}{1-q} \Big)^{z-1} f_1 \Big( \frac{q^{k+1}}{1-q} \Big) \right)\left( \sum_{n=0}^\infty q^n \Big( \frac{q^n}{1-q} \Big)^{z-1} f_2 \Big( \frac{q^{n+1}}{1-q} \Big) \right)\\
= \sum_{k, n=0}^\infty q^{k+n} \Big( \frac{q^{k+n}}{(1-q)^2} \Big)^{z-1} f_1 \Big( \frac{q^{k+1}}{1-q} \Big)   f_2 \Big( \frac{q^{n+1}}{1-q} \Big) \\
= \Big( \frac{1}{1-q} \Big)^{z-1} \sum_{m=0}^\infty q^{m} \Big( \frac{q^{m}}{1-q} \Big)^{z-1} \underbrace{\left( \sum_{k=0}^m f_1 \Big( \frac{q^{k+1}}{1-q} \Big)   f_2 \Big( \frac{q^{m+1-k}}{1-q} \Big) \right)}_{:= f_1 \circ f_2 \left( q\frac{q^m}{1-q} \right)}.
\end{gather*}
\end{proof}

~





For the multiplication operator $M_z$ we have the following fact.
\begin{proposition}
$M_z$ is bounded from $\mathbf H_{2,q}$ into itself with norm $\|M_z\|\le \frac{1}{1-q}$.
\end{proposition}

\begin{proof}
  This follows from
  \[
\frac{\frac{1}{1-q}-z\overline{w}}{\prod_{j=0}^\infty (1-z\overline{w}(1-q)q^j)}=\frac{1}{1-q}\frac{1}{\prod_{j=1}^\infty (1-z\overline{w}(1-q)q^j)}.
\]
Since the kernel $\frac{1}{1-q}\frac{1}{\prod_{j=1}^\infty (1-z\overline{w}(1-q)q^j)}$ is positive definite in $\mathbb D_{1/1-q}$ so is the kernel
\[
  \frac{\frac{1}{1-q}-z\overline{w}}{\prod_{j=0}^\infty (1-z\overline{w}(1-q)q^j)},
\]
and we conclude with the characterization of multipliers in a reproducing kernel Hilbert space.
  \end{proof}

\begin{lemma} (see e.g. \cite[Exercise 4.2.25, pp. 165 and 185]{MR3560222})
  \begin{equation}
    (R_qf)(z)=\lambda f(z)\quad\iff\quad f(z)=\frac{c}{\prod_{j=0}^\infty(1-\lambda (1-q)zq^j)}.
    \end{equation}
  \end{lemma}

\begin{proposition} The $q$-exponential satisfy
  \begin{equation}
    \label{voltaire}
    (    R_qE_q(\cdot \overline{w}))(z)=\overline{w}E_q(z\overline{w}).
    \end{equation}
\end{proposition}

\begin{proof}
  We note that $
E_q(qz\overline{w})=(1-z\overline{w}(1-q))E_q(z\overline{w})$ and so
\[
  \begin{split}
    (    R_qE_q(\cdot \overline{w}))(z)&=\frac{E_q(z\overline{w})-E_q(qz\overline{w})}{(1-q)z}\\
    &=\frac{E_q(z\overline{w})-(1-z\overline{w}(1-q))E_q(z\overline{w})}{(1-q)z}\\
    &=\overline{w}E_q(z\overline{w}).
\end{split}
  \]
  \end{proof}

\begin{theorem} Let $q\in[0,1)$. The only Hilbert space of functions which is analytic in a neighborhood of the origin and for which
\begin{equation}
R_q^*=M_z
\end{equation}
is $\mathbf H_{2,q}$ (up to a multiplicative factor for the inner product).
\end{theorem}
\begin{proof}
We have that $E_q(z \overline w) = K_q(z, w)$ is the reproducing for $\mathbf{H}_{2,q},$ i.e. $f(z) = \langle f, E_q(\cdot\overline{z})\rangle_{\mathbf H_{2,q}}$. Using \eqref{voltaire} we can write:
  \[
    \begin{split}
      \left(R_q^*E_q(\cdot\overline{w})\right)(z)&=      \langle R_q^*E_q(\cdot\overline{w}),E_q(\cdot\overline{z})\rangle_{\mathbf H_{2,q}}\\
      &=      \langle E_q(\cdot\overline{w}), R_qE_q(\cdot\overline{z})\rangle_{\mathbf H_{2,q}}\\
      &=       \langle E_q(\cdot\overline{w}), \overline{z}E_q(\cdot\overline{z})\rangle_{\mathbf H_{2,q}}\\
      &=z E_q(z\overline{w}).
    \end{split}
    \]
  \end{proof}

\begin{proposition}\label{Prop3.9}
The space $\mathbf H_{2,q}$ is a de Branges-Rovnyak space.
\end{proposition}

\begin{proof}
  This follows from \cite[Theorem 2.1 p. 51]{MR4032209}, since
  the sequence $\fact{k}$, for $k=0,1,\ldots$ is an increasing sequence with initial term $1$.
\end{proof}

We now compute the adjoint of $R_q$ in $\mathbf H_{2,q}$. Since
\begin{equation}
\inner{z^n, z^m}_{\mathbf H_{2,q}}= \fact{n} \, \delta_{n, m},
  \end{equation} we have
\begin{gather*}
\inner{z^n, R_q z^m}_{\mathbf H_{2,q}} = \inner{z^n, \frac{z^m-q^m z^m}{(1-q)z}}_{\mathbf H_{2,q}} = (1+q+\cdots +q^{m-1}) \inner{z^n, z^{m-1}}_{\mathbf H_{2,q}} \\
= (1+q+\cdots +q^{n})\fact{n} = \fact{n+1}\\
 =  \inner{z^{n+1}, z^m}_{\mathbf H_{2,q}} =: \inner{R^\ast_q z^n, z^m}_{\mathbf H_{2,q}}.
\end{gather*}
Therefore, we obtain  $R_q^*=M_z.$ \\

In the case $q=1$ the Fock space can be characterized (up to a multiplicative positive factor in the inner product) as the only Hilbert space of power series converging in a convex neighborhood of the origin and such that
\begin{equation}
  (R_0^*f)(z)=\int_{[0,z]}f(s)ds,
  \label{R0*int}
\end{equation}that is, $R_0^\ast$ coincides with the integration operator. It is therefore natural to try and define the integral in $\mathbf H_{2,q}$ by $R_q^*$ for $q\in(0,1)$.  

\begin{lemma}
  The operator $R_0$ is bounded in $\mathbf H_{2,q}$ and it holds that (with $e_k(z)=z^k$)
  \begin{equation}
R_0^*e_k=\frac{e_{k+1}}{1+q+\cdots +q^k},\quad k=0,1,\ldots
    \end{equation}
  \end{lemma}

 \begin{proof}
   We have for $k\ge 1$ and $\ell\ge 0$
   \[
     \begin{split}
       \langle R_0e_{k},e_{\ell}\rangle_{\mathbf H_{2,q}}&=\langle e_{k-1},e_{\ell}\rangle_{\mathbf H_{2,q}}\\
       &=\delta_{k-1,\ell} \fact{\ell} \\
       &=\delta_{k-1,\ell}\langle e_k,e_k\rangle_{\mathbf H_{2,q}}\frac{\fact{\ell}}{\fact{k}}     \\ 
              &=\delta_{k-1,\ell}\langle e_k,e_k\rangle_{\mathbf H_{2,q}}\frac{1}{1+q+\cdots+q^{\ell}}\\
&=\langle e_k ,R_0^*e_\ell\rangle_{\mathbf H_{2,q}},
     \end{split}
   \]
   with
   \begin{equation}
     \label{R0star}
R_0^*e_\ell=\frac{e_{\ell+1}}{1+q+\cdots+q^{\ell}}.
     \end{equation}
\end{proof}

Consider the $q$-Jackson integral
$$
\int_0^a f(x)d_qx:=(1-q)a\sum_{k=0}^\infty q^k f(q^k a),
$$
which is said to converge provided that the sum on the right-hand-side converges absolutely.

\begin{lemma}
  \label{lemma-int}
  \begin{equation}
    \label{new345}
\int_0^{z}x^\ell d_qx=z^{\ell+1}\frac{1}{1+q+\cdots +q^\ell}.
    \end{equation}
  \end{lemma}

  \begin{proof}
By definition we have
\[
  \begin{split}
    \int_0^{z}x^\ell d_qx&=z(1-q)\sum_{k=0}^\infty q^k(q^kz)^\ell =z^{\ell+1}(1-q)\left(\sum_{k=0}^\infty (q^{1+\ell})^k\right)\\
    &=z^{\ell+1}(1-q)\frac{1}{1-q^{\ell+1}} =z^{\ell+1}\frac{1}{1+q+\cdots +q^\ell}.
    \end{split}
  \] \end{proof}

It is well know that
\begin{equation}
\label{ertyuiop}
\partial^*=M_z
\end{equation}
in the Fock space, and that in fact the Fock space is characterized (up to a positive multiplicative constant in the inner product)
by this equality; see \cite{MR0157250}. In \cite{nati_01} it is proved that in the Hardy space we have
\begin{equation}
  \partial^*=M_z\partial M_z,
\label{nati001}
  \end{equation}
  and that the above equality does characterize the Hardy space (as usual, up to a positive multiplicative constant in the inner product). We now prove a formula which
  is valid for $q\in[0,1]$ and englobes the two above formulas.

  \begin{theorem}
    Let $q\in[0,1]$.
    Then in $\mathbf H_{2,q}$ it holds that
    \begin{equation}
      \label{strange123}
\partial^*=M_z\partial R_0^*
\end{equation}
and this equality characterizes the space $\mathbf H_{2,q}$ up to a positive multiplicative constant in the inner product.
\label{structure}
\end{theorem}

When $q=0$ (Hardy space) we have $R_0^* $ that $R_0^*=M_z$ and so \eqref{strange123} reduces to
\[
M_z\partial M_z,
\]
i.e. \eqref{nati001}. When $q=1$ (Fock space), we have $R_0^*=I$ (the integration operator) and $\partial I e_k=e_k$, $k=0,1,\ldots$. We thus get back \eqref{ertyuiop}.

    \begin{proof}[Proof of Theorem \ref{structure}]
Let $k\in\mathbb N_0$. Let us set {\it a priori} $\partial^* e_k=a_{k,q} e_{k+1}$  for some $a_{k,q}\in\mathbb C$. We have on the one hand
\[
    \begin{split}
          \langle \partial^*e_k, e_{k+1}\rangle_{\mathbf H_{2,q}}&=\langle e_k, \partial e_{k+1}\rangle_{\mathbf H_{2,q}}\\
          &=(k+1)\langle e_k,  e_{k}\rangle_{\mathbf H_{2,q}}\\
          &=(k+1) \fact{k} 
        \end{split}
  \]
 and on the other hand, with $\partial^*e_k=a_{k,q}e_{k+1}$ we have
\[
        \begin{split}
          \langle \partial^*e_k, e_{k+1}\rangle_{\mathbf H_{2,q}}
          &=a_{k,q}\langle e_{k+1},  e_{k+1}\rangle_{\mathbf H_{2,q}}\\
          &=a_{k,q} \fact{k+1}.
        \end{split}
        \]
        Thus
        \[
a_{k,q}\fact{k+1} =          (k+1) \fact{k}
\]
from which we get
\begin{equation}
  \label{la-rochelle}
  a_{k,q}=\frac{k+1}{1+q+\cdots +q^k}.
\end{equation}
In view of \eqref{R0star}, we can write
\[
  \begin{split}
    \partial^*e_k&=\frac{(k+1)M_ze_k} {1+q+\cdots +q^k}\\
    &=(k+1)R_0^*e_k\\
    &=M_z\partial R_0^*e_k
    \end{split}
  \]
  since
\[
  M_z \partial e_{k+1}=(k+1)e_{k+1},\quad k=0,1,\ldots
  \]
      \end{proof}

Note that in \eqref{la-rochelle}, we set $a_{k,0}=k+1$  and $a_{k,1}=1,$      as it should be.

\begin{theorem}We have
  \begin{equation}
    \label{mzqrq}
    M_z^*=R_qM_zR_0.
    \end{equation}
\end{theorem}

\begin{proof} For $m=1,2, \ldots$  we get
\begin{gather*}
\inner{z^n, R_qM_zR_0 z^m}_{\mathbf H_{2,q}} = \inner{z^n, R_qM_z z^{m-1}}_{\mathbf H_{2,q}} = \inner{z^n, R_q z^{m}}_{\mathbf H_{2,q}}.
\end{gather*} By Proposition \ref{Prop3.9} we obtain
\begin{gather*}
\inner{z^n, R_qM_zR_0 z^m}_{\mathbf H_{2,q}} = \inner{z^n, R_q z^{m}}_{\mathbf H_{2,q}} = \inner{M_z z^n,  z^m}_{\mathbf H_{2,q}}.
\end{gather*} We conclude our proof with the observation that $0=\inner{z^n, R_qM_zR_0 z^0}_{\mathbf H_{2,q}} = \inner{M_z z^n,  z^0}_{\mathbf H_{2,q}}.$
  \end{proof}

\section{The space ${\mathcal F_{2,q}}$}
\setcounter{equation}{0}

The space $\mathcal F_{2,q}$ appeared in \cite{daf1} motivated by a study of discrete analytic functions.
\begin{definition}
Consider the reproducing kernel
$$
K_{2,q}(z,w)=\sum_{n=0}^\infty \frac{z^n\overline{w}^n}{([n_q]!)^2}.
$$
Then the corresponding reproducing kernel Hilbert space $\mathcal{F}_{2,q}$ is the space of all functions $f(z) =\sum_{n=0}^\infty f_n z^n$ such that
$\sum_{n=0}^\infty |f_n|^2 ([n_q]!)^2<\infty$.
\end{definition}

In this way, we have $K_q =: K_{1, q}$ and $K_{2,q}$ as the reproducing kernels of $\mathbf{H}_{2,q}$ and $\mathcal{F}_{2,q},$ respectively. As both kernels are positive definite and the same holds for its difference $K_{1, q}-K_{2, q}$ we get that ${\mathcal F_{2,q}}$ is contractively included in ${\mathbf H_{2,q}}$ (see  \cite{MR2002b:47144, Aronszajn44}).

\begin{remark}
We observe that for $f_1(z) = f_2(z) =  E_q^{-1}(z)$ we have
\begin{eqnarray*}
\left[ \mathcal M_q(E_q^{-1})(n+1) \right]^2  &= &
\Big( \frac{1}{1-q} \Big)^{n} \sum_{m=0}^\infty q^{m} \Big( \frac{q^{m}}{1-q} \Big)^{z-1} \left[ \sum_{k=0}^m  E_q^{-1} \Big( \frac{q^{k+1}}{1-q} \Big)   E_q^{-1} \Big( \frac{q^{m+1-k}}{1-q} \Big) \right] \\
(\fact{n})^2 &= & \Big( \frac{1}{1-q} \Big)^{n} \sum_{m=0}^\infty q^{m} \Big( \frac{q^{m}}{1-q} \Big)^{n} \sum_{k=0}^m (q^{k+1}; q)_\infty  (q^{m+1-k}; q)_\infty.
\end{eqnarray*}

Hence, we get as density $\omega_{2,q}$ of our $q$-Fock space $\mathcal F_{2,q}$
\begin{equation} \label{Eq:3.0002}
\omega_{2,q}(|z|) := \left(\frac{1}{1-q} \right)^{|z|^2-1} (E_q^{-1} \circ E_q^{-1})(|z|^2),
\end{equation} and satisfying to
\begin{equation} \label{Eq:3.0003}
\mathcal M_q (\omega_{2,q})(n+1) =  \left(\frac{1}{1-q} \right)^{n} \mathcal M_q (E_q^{-1} \circ E_q^{-1})(n+1) = (\fact{n})^2.
\end{equation}
\end{remark}

Now, we can define $T_q : \mathbf H_2 \mapsto \mathcal F_{2,q}$ given as $z^n \to \frac{z^n}{\fact{n}}.$

\begin{lemma}\label{Lem:4.2} In $\mathbf H_2$ it holds:
\begin{equation}
R_q T_q =T_qR_0.
  \end{equation}
\end{lemma}
\begin{proof} The case of $n=0$ is immediate. For $n=1, 2, \ldots$ we have
\begin{gather*}
R_q T_q z^n = R_q \left( \frac{z^n}{\fact{n}}  \right) = \frac{1}{\fact{n}} (1+q+\cdots +q^{n-1}) z^{n-1}\\
= \frac{1}{\fact{n-1}}z^{n-1} = T_q z^{n-1} = T_q R_0 z^{n}.
\end{gather*}
\end{proof}

\begin{theorem}\label{Th:4.3}
The map $T_q$ is an isometry from $\mathbf H_{2}$ onto ${\mathcal  F_{2,q}}$.
  \end{theorem}

 \begin{proof} We have $\langle e_n, e_m\rangle_{\mathcal  F_{2,q}}= (\fact{n})^2  \delta_{n,m}.$ Hence, we get
   \[
     \begin{split}
       \langle T_qe_n,T_qe_m\rangle_{\mathcal  F_{2,q}}&=\frac{1}{\fact{n} }\frac{1}{\fact{m}}\langle e_n,e_m\rangle_{\mathcal  F_{2,q}}\\
       &=\delta_{n,m}\frac{(\fact{n})^2}{(\fact{n})^2}\\
       &=\delta_{n,m}\\
       &=\langle e_n,e_m\rangle_{\mathbf H_{2}}.
     \end{split}
     \]
\end{proof}

\begin{theorem} \label{Th:4.4}
  In $\mathcal F_{2,q}$ it holds that
  \begin{equation}
\label{Rqstar}
    R_q^*e_n=\frac{e_{n+1}}{[n]_q},\quad n=0,1,\ldots
  \end{equation}
    and
    \begin{equation}
      \label{paris}
\mathcal{I}-R_q^*R_q=C^*C,
\end{equation}
and this structural identity characterizes the space ${\mathcal F_{2,q}}$ up to a multiplicative factor.
  \end{theorem}

  \begin{proof}
    To prove \eqref{Rqstar} we write
\[
    \begin{split}
      \langle R_q^*e_n,e_m\rangle_{\mathcal F_{q,2}}&=\langle e_n,R_qe_m\rangle_{\mathcal F_{2,q}}\\
      &=[m]_q\langle e_n,e_{m-1}\rangle_{\mathcal F_{q,2}}\\
      &=\delta_{m-1,n}([n]_q!)^2[m]_q.
    \end{split}
  \]
  On the other hand we show that one can assume that $R_q^*e_n=\alpha_n e_{n+1}$; we have
  \[
    \begin{split}
      \langle R_q^*e_n,e_m\rangle&=\alpha_n\langle e_{n+1},e_m\rangle_{\mathcal F_{2,q}}\\
      &=\alpha_n\delta_{n+1,m}([n+1]_q!)^2.
    \end{split}
  \]
  Comparing these equalities we obtain
  \[
\alpha_n([n+1]_q!)^2=([n]_q!)^2[n]_q,
\]
so that $\alpha_n=\frac{1}{[n]_q}$. It follows that
\[
  R_q^*R_qe_n=\begin{cases}
\, 0,\,\hspace{.9mm}\quad n=0,\\
    \, e_n,\quad n=1,2,\ldots
    \end{cases}
  \]
  and hence the result.
    \end{proof}

From the previous computations we also have:

\begin{proposition}
$R_0^*$ is an isometry in ${\mathcal F_{2,q}}$.
\end{proposition}

\begin{proof} This is a direct consequence of the fact that
\[
  R_0R_0^*e_n=e_n,
\]for all $n \in \mathbb N_0.$
  \end{proof}

From Lemma \ref{lemma-int} we have:

\begin{proposition}
  In $\mathcal F_{2,q},$ it holds
  \[
R_q^*=\mathrm I,
\]
where $\mathrm I$ is the integration operator.
  \end{proposition}

 We now use well a known method in characteristic function theory (see e.g. \cite{adrs} in the case of Pontryagin spaces) and rewrite \eqref{paris} as
  \[
\begin{pmatrix}R_q\\ C\end{pmatrix}^*\begin{pmatrix}R_q\\ C\end{pmatrix}= \mathcal I.
\]
The operator
\[
\begin{pmatrix}\mathcal I & 0\\ 0& 1\end{pmatrix}
- {\begin{pmatrix}R_q\\ C\end{pmatrix}\begin{pmatrix}R_q\\ C\end{pmatrix} }^*
\]
is therefore positive and for instance using its square root, one can find a Hilbert space $\tilde{\mathcal H}$ and operators $B$ and $D,$

\[
\begin{pmatrix}B\\ D\end{pmatrix}\,\,:\,\, \tilde{\mathcal H} \longrightarrow \mathcal F_{2,q}\oplus\mathbb C,
  \]
such that
    \[
\begin{pmatrix}\mathcal I  & 0\\ 0& 1\end{pmatrix}   -\begin{pmatrix}R_q\\ C\end{pmatrix}\begin{pmatrix}R_q\\ C\end{pmatrix}^*=\begin{pmatrix}B\\ D\end{pmatrix}\begin{pmatrix}B\\ D\end{pmatrix}^*.
      \]
      The operator matrix
      \begin{equation}
        \label{abcd}
        \begin{pmatrix}R_q&B\\
          C&D\end{pmatrix}
      \end{equation}
      is co-isometric. We set
\begin{equation}\label{Sq}
S_q(z)=D+z C (\mathcal I-zR_q)^{-1}B.
\end{equation} 

We now look into the properties of the matrix (\ref{abcd}). We observe that
$$ \begin{pmatrix}
	 \mathcal I & 0\\
          0& 1 \end{pmatrix}=  \begin{pmatrix}R_q&B\\
          C&D\end{pmatrix}
        \begin{pmatrix}R_q&B\\
          C&D\end{pmatrix}^\ast =  \begin{pmatrix}R_q R_q^\ast + B B^\ast & R_q C^\ast + B D^\ast \\
          CR^\ast_q +D B^\ast &CC^\ast + D D^\ast \end{pmatrix},$$
so that we get $$DD^\ast = 1 - CC^\ast, \quad B B^\ast = \mathcal I  - R_q R_q^\ast, \quad  B D^\ast = - R_q C^\ast.$$ Hence, from (\ref{Sq}) we get
\begin{gather*}
S_q(z) [S_q(w)]^\ast=[D\mathcal I+z C (\mathcal I-zR_q)^{-1}B] [D^\ast \mathcal I+ \ov w  B^\ast  (\mathcal I-\ov w R_q^\ast)^{-1} C^\ast ]  \\
= D D^\ast \mathcal I  + z C (\mathcal I-zR_q)^{-1}BD^\ast + \ov w  D B^\ast  (\mathcal I-\ov w R_q^\ast)^{-1} C^\ast  + z \ov w C (\mathcal I-zR_q)^{-1}B B^\ast  (\mathcal I-\ov w R_q^\ast)^{-1} C^\ast  \\
= (1 - CC^\ast)\mathcal I   + z C (\mathcal I-zR_q)^{-1}BD^\ast + \ov w  D B^\ast  (\mathcal I-\ov w R_q^\ast)^{-1} C^\ast  + z \ov w C (\mathcal I-zR_q)^{-1}B B^\ast  (\mathcal I-\ov w R_q^\ast)^{-1} C^\ast
\end{gather*}so that
\begin{gather*}
\mathcal I -S_q(z) [S_q(w)]^\ast \\
=  CC^\ast \mathcal I  - z C (\mathcal I-zR_q)^{-1}BD^\ast - \ov w  D B^\ast  (\mathcal I-\ov w R_q^\ast)^{-1} C^\ast  - z \ov w C (\mathcal I-zR_q)^{-1}B B^\ast  (\mathcal I-\ov w R_q^\ast)^{-1} C^\ast \\
=  CC^\ast \mathcal I + z C (\mathcal I-zR_q)^{-1}R_q C^\ast  + \ov w  C R_q^\ast (\mathcal I-\ov w R_q^\ast)^{-1} C^\ast   - z \ov w C (\mathcal I-zR_q)^{-1} (\mathcal I - R_q R_q^\ast )  (\mathcal I-\ov w R_q^\ast)^{-1} C^\ast \\
=  C(\mathcal I-zR_q)^{-1}  \Big[ \underbrace{(\mathcal I-zR_q)(1-\ov w R_q^\ast)  + z  R_q (\mathcal I-\ov w R_q^\ast) + \ov w (\mathcal I-zR_q)  R_q^\ast -  z \ov w  (\mathcal I - R_q R_q^\ast )}_{(A)} \Big] (\mathcal I-\ov w R_q^\ast)^{-1}  C^\ast.
\end{gather*} Easy calculations give now
\begin{eqnarray*}
(A) &=&  (\mathcal I-zR_q)(\mathcal I-\ov w R_q^\ast)  + z  R_q (\mathcal I-\ov w R_q^\ast) + \ov w (\mathcal I-zR_q)  R_q^\ast -  z \ov w  (\mathcal I - R_q R_q^\ast ) \\
&=& \mathcal I - zR_q -\ov w R_q^\ast + z \ov w R_qR_q^\ast + z  R_q - z \ov w   R_q  R_q^\ast + \ov w R_q^\ast - z \ov w R_q R_q^\ast -  z \ov w \mathcal I  + z \ov w R_q R_q^\ast \\
&=&(1 - z \ov w) \mathcal I,
\end{eqnarray*}

Hence, it holds that
        \begin{equation}
          \label{Schur_q}
          \frac{\mathcal I-S_q(z)S_q(w)^*}{1-z\overline{w}}=C(\mathcal I-zR_q)^{-1}[(\mathcal I-wR_q)^*]^{-1}C^*,\quad z,w\in\mathbb D.
        \end{equation}

 The operator $S_q$ bears various names in operator theory; it is the characteristic operator function, or the transfer function, or the scattering function, associated to the operator matrix \eqref{abcd}. From \eqref{Schur_q} one sees that $S_q$ is analytic and contractive in the open unit
 disk, i.e. is a Schur function.\\

        When $q=0$ we have for $ f(z)=\sum_{n=0}^\infty c_nz^n$ that
        \[
          C(\mathcal I-zR_0)^{-1}f =f(z),\quad z \in\mathbb D.
          \]
Here, for $0<q \leq 1$ we define
        \[
f_q(z)=C(\mathcal I-zR_q)^{-1}f,\quad z \in\mathbb D.
          \]
As $CR_q^nf=[n]_q! \, c_n$ we get that the coefficients
        \begin{equation} \label{Eq:4.001}
c_n=\frac{CR_q^nf}{[n]_q!}
\end{equation}
are independent of $q$ and one has $f_q(z)=f(z),$ that is, we obtain $f(z) = C(\mathcal I-zR_q)^{-1}f$  for all $z \in \mathbb D$.\\

\section{Aknowledgements}
D. Alpay thanks the Foster G. and Mary McGaw Professorship in Mathematical Sciences, which supported his research. The second and third author were supported by CIDMA, through the Portuguese FCT (UIDP/04106/2020 and UIDB/04106/2020). 
  
The authors thank warmly Jeanne Scott and Martin Schork for pointing out references~\cite{Milne1978} and~\cite{Schork2016}, respectively, where the $q$-Stirling numbers appearing in this paper can already be found.

%

      \bibliographystyle{plain}
\def\cprime{$'$} \def\cprime{$'$} \def\cprime{$'$}
  \def\lfhook#1{\setbox0=\hbox{#1}{\ooalign{\hidewidth
  \lower1.5ex\hbox{'}\hidewidth\crcr\unhbox0}}} \def\cprime{$'$}
  \def\cprime{$'$} \def\cprime{$'$} \def\cprime{$'$} \def\cprime{$'$}
  \def\cprime{$'$}


\begin{thebibliography}{10}

\bibitem{MR2002b:47144}
D.~Alpay.
\newblock {\em The {S}chur algorithm, reproducing kernel spaces and system
  theory}.
\newblock American Mathematical Society, Providence, RI, 2001.
\newblock Translated from the 1998 French original by Stephen S. Wilson,
  Panoramas et Synth\`eses.

\bibitem{MR3560222}
D.~Alpay.
\newblock {\em A complex analysis problem book}.
\newblock Birkh\"auser/Springer, Cham, 2016.
\newblock Second edition.

\bibitem{Trevor23}
D.~Alpay, P.~Cerejeiras, U.~Kaehler, and T.~Kling.
\newblock Commutators on {F}ock spaces.
\newblock {\em J. {M}ath. {P}hys.}, 64:042102--21pages, 2023.

\bibitem{MR4032209}
D.~Alpay, F.~Colombo, and I.~Sabadini.
\newblock The {F}ock space as a {D}e {B}ranges-{R}ovnyak space.
\newblock {\em Integral Equations Operator Theory}, 91(6):Paper No. 51, 12,
  2019.

\bibitem{adrs}
D.~Alpay, A.~Dijksma, J.~Rovnyak, and H.~de~Snoo.
\newblock {\em {Schur} functions, operator colligations, and reproducing kernel
  {P}ontryagin spaces}, volume~96 of {\em Operator theory: {A}dvances and
  {A}pplications}.
\newblock Birkh{\" a}user Verlag, Basel, 1997.

\bibitem{daf1}
D.~{Alpay}, P.~{Jorgensen}, R.~{Seager}, and D.~{Volok}.
\newblock {On discrete analytic functions: Products, rational functions and
  reproducing kernels}.
\newblock {\em Journal of Applied Mathematics and Computing}, 41:393--426,
  2013.

\bibitem{MR3816055}
D.~Alpay and M.~Porat.
\newblock Generalized {F}ock spaces and the {S}tirling numbers.
\newblock {\em J. Math. Phys.}, 59(6):063509, 12, 2018.

\bibitem{nati_01}
N.~Alpay.
\newblock A new characterization of the {H}ardy space and of other spaces of
  analytic functions.
\newblock {\em ArXiv (2020); To appear in İstanb. Univ., Sci. Fac., J. Math.
  Phys. Astron.}, pages 1--10, 2023.

\bibitem{Aronszajn44}
N.~Aronszajn.
\newblock La th{\'e}orie g{\'e}n{\'e}rale des noyaux reproduisants et ses
  applications.
\newblock {\em Math. Proc. Cambridge Phil. Soc.}, 39:133--153, 1944.

\bibitem{MR0157250}
V.~Bargmann.
\newblock On a {H}ilbert space of analytic functions and an associated integral
  transform.
\newblock {\em Comm. Pure Appl. Math.}, 14:187--214, 1961.

\bibitem{bargmann}
V.~Bargmann.
\newblock Remarks on a {H}ilbert space of analytic functions.
\newblock {\em Proceedings of the {N}ational {A}cademy of {Arts}}, 48:199--204,
  1962.

\bibitem{Cai2}
Y.~Cai, R.~Ehrenborg, and M.~Readdy.
\newblock $q$-{S}tirling identities revisited.
\newblock {\em {E}lectron. {J}. {C}omb.}, 25:1--37, 2018.

\bibitem{Cai1}
Y.~Cai and M.~Readdy.
\newblock $q$-{S}tirling numbers: a new view.
\newblock {\em {A}dv. {A}ppl. {M}ath.}, 86:50--80, 2017.

\bibitem{duren}
P.L. Duren.
\newblock {\em Theory of $H^p$ spaces}.
\newblock Academic press, New York, 1970.

\bibitem{MR1102893}
K.~Hoffman.
\newblock {\em Banach spaces of analytic functions}.
\newblock Dover Publications Inc., New York, 1988.
\newblock Reprint of the 1962 original.

\bibitem{Jackson1910}
F.~H. Jackson.
\newblock On $q$-definite integrals.
\newblock {\em Quart. J.}, 41:193--203, 1910.

\bibitem{Kac2001}
V.~G. Kac and P.~Cheung.
\newblock {\em {Q}uantum calculus}.
\newblock Universitext (UTX). Springer, New York, 2001.

\bibitem{Schork2016}
T.~Mansour and M.~Schork.
\newblock {\em Commutation Relations, Normal Ordering, and Stirling Numbers}.
\newblock Chapman and Hall/CRC, New York, 2016.

\bibitem{Milne1978}
S.C. Milne.
\newblock A q-analog of restricted growth functions, {D}obinkski's equality,
  and {C}harlier polynomials.
\newblock {\em Trans. Amer. Math. Soc.}, 245:89--118, 1978.

\bibitem{MR924157}
W.~Rudin.
\newblock {\em Real and complex analysis}.
\newblock McGraw-Hill Book Co., New York, third edition, 1987.

\bibitem{Leeuwen1995}
H.~van Leeuwen and H.~Maassen.
\newblock A $q$-deformation of the {G}auss distribution.
\newblock {\em J. {M}ath. {P}hys.}, 36:4743--4756, 1995.

\bibitem{Zhu2012}
K.~Zhu.
\newblock {\em {A}nalysis on {F}ock spaces}, volume 263 of {\em Graduate Texts
  in Mathematics (GTM)}.
\newblock Springer, New York, NY, 2012.

\end{thebibliography}

\end{document}